\documentclass[reqno,11pt]{amsart}
\usepackage[top=2cm,bottom=2.5cm,right=2.5cm,left=2.5cm]{geometry}
\usepackage[english]{babel}
\usepackage[T1]{fontenc}
\usepackage{times}
\usepackage{mathtools, bm}
\usepackage{amssymb, bm}
\usepackage{graphicx}
\usepackage{mathrsfs}
\usepackage{stmaryrd}
\usepackage{amsthm}
\usepackage{listings}
\usepackage{pict2e}
\usepackage{pgf, tikz}
\usepackage{dsfont}
\usepackage{algorithm2e}
\usepackage{algorithmic}
\usepackage{multicol}
\usepackage{hyperref}
\usepackage{bbm}

\newtheorem{theorem}{Theorem}[section]

\newtheorem{proposition}[theorem]{Proposition}



\author{Jules~Grass}
\address{Université Claude Bernard Lyon 1, CNRS UMR 5208, Institut Camille Jordan, 69622 Villeurbanne, France. \url{grass@math.univ-lyon1.fr}}

\author{Arnaud~Guillin}
\address{Laboratoire de Mathématiques Blaise Pascal, CNRS UMR 6620, Université Clermont-Auvergne, avenue des Landais, F-63177 Aubière. \url{arnaud.guillin@uca.fr}}

\author{Christophe~Poquet}
\address{Université Claude Bernard Lyon 1, CNRS UMR 5208, Institut Camille Jordan, 69622 Villeurbanne, France. \url{poquet@math.univ-lyon1.fr}}

\title[Sharp propagation of chaos with non constant diffusion coefficient]{Sharp propagation of chaos for McKean-Vlasov equation with non constant diffusion coefficient}

\newcommand{\fk}{\mu_{t}^{k,N}}
\newcommand{\ftk}{\mu_{t}^{\otimes k}}

\begin{document}

\begin{abstract}
We present a method to obtain sharp local propagation of chaos results for a system of N particles with a diffusion coefficient that it not constant and may depend of the empirical measure. This extends the recent works of Lacker \cite{Lacker} and Wang \cite{Wang} to the case of non constant diffusions. The proof relies on the BBGKY hierarchy to obtain a system of differential inequalities on the relative entropy of k particles, involving the fisher information.
\end{abstract}

\maketitle

\section{Introduction and main result}
In this work, we are interested in the following system of N particles on the $d$-dimensional torus $\mathbb{T}^{d}$
\begin{equation} 
dV^{i,N}_{t}=\frac{1}{N} \sum_{j=1}^{N} b(V^{i,N}_{t}-V^{j,N}_{t})dt+\sqrt{2} \Bigl(a_{1}(V_{t}^{i,N})+\frac{1}{N} \sum_{j=1}^{N} a_{2}(V_{t}^{i,N}-V_{t}^{j,N}) \Bigl)^{\frac{1}{2}} dB^{i,N}_{t}
\label{(1)},
\end{equation}
where $(B^{1,N},...,B^{N,N})$ are independent Brownian motions. We suppose that the matrix valued functions $a_{1}$ and $a_{2}$ are symmetric and satisfy
\begin{equation} \forall (v_{1},v_{2}) \in ({T}^{d})^{2},\,  a_{1}(v_{1})+a_{2}(v_{1}-v_{2}) \geq \lambda_{1} \text{Id}.\label{condition a} \end{equation}
We also assume that $b$, $\nabla \cdot a_{1}$, $a_{2}$ and $\nabla \cdot a_{2}$ are bounded. Throughout this paper, we suppose moreover that the N particles in the dynamics \ref{(1)} are exchangeable, that is, for all permutation $\sigma$ of $\llbracket 1,N \rrbracket$, we have $\text{Law}\bigl(V_{t}^{\sigma(1),N},...,V_{t}^{\sigma(N),N}\bigl)=\text{Law}\bigl(V_{t}^{1,N},...,V_{t}^{N,N}\bigl)$. 

An example of particular interest of such a system is when one takes $b:=\nabla \cdot a_{2}$ and $a_{1}=0$, obtaining a Landau-like equation that has notably been recently studied by Carillo, Guo and Jabin \cite{Equation}.

We are interested in the behavior of (\ref{(1)}) when $N \to \infty$ and especially the property known as propagation of chaos, i.e, the convergence of $\fk:=\text{Law}\bigl(V_{t}^{1,N},...,V_{t}^{k,N} \bigl)$ to $\ftk$, where $\mu_t$ denotes the law of the solution to the McKean-Vlasov equation
\begin{equation}
\left\{
            \begin{array}{ll}
                dV_{t}= b \ast \mu_{t}(V_{t})dt+\sqrt{2}\bigl(a_{1}(V_{t})+a_{2} \ast \mu_{t}(V_{t})\bigl)^{\frac{1}{2}} dB_{t} \\ \mu_{t}=\text{Law}(V_{t})
            \end{array}
            \right. .
\label{(2)}
\end{equation}
This property has garnered a lot of attention from the mathematical community ever since the second half of the 20th century, although its root can be traced back to the birth of statistical mechanics with the assumption of molecular chaos (also know as the Stosszahlansatz). Several methods have been used troughout the years to prove propagation of chaos, starting with compactness arguments \cite{Ref compacite 1, Ref compacite 2, M96}, and then coupling arguments, see \cite{Ref couplage 1,Ref couplage 2,E16} for historical references or to \cite{M03,CGM08,DEGZ20} for uniform in time results in the case of convex or non non convex interactions and \cite{GLBM23} for the 1D Coulomb case. Recently, entropy methods have gained a lot of attraction, notably because they were successfully applied to systems with singular interaction \cite{Jabin wang}. For a review of those methods see \cite{Review POC 1,Review POC 2}. With the coupling or entropy approaches, one typically obtains a convergence rate of $\fk$ to $\ftk$ of order $\mathcal{O} \bigl( \sqrt{\frac{k}{N}} \bigl)$ in total variation or Wasserstein distance. For the relative entropy method this is done by proving that $H(\mu_{t}^{N,N}|\mu_{t}^{\otimes N})=\mathcal{O}(1)$, and then concluding by using the subadditivity of relative entropy. Let us also mention the recent modulated energy approach for the singular case in \cite{RS23}, and further successfully developped for attractive interactions \cite{BJW23}, however with non comparable rates.

Recently, under some assumptions on the interaction term, that are for example valid if the interaction is bounded or Lipschitz continuous, Lacker \cite{Lacker} was able to obtain optimal bounds for $H_{t}^{k}:=H(\fk|\ftk)$, of order $\mathcal{O}\bigl( \frac{k^{2}}{N^{2}} \bigl)$ (using Pinsker's inequality this means a convergence rate of $\mathcal{O}\bigl( \frac{k}{N} \bigl)$ for the total variation distance). He proved moreover that such a bound is optimal for a simple Gaussian system. The novelty of this approach is its local character, by estimating $H_{t}^{k}$ directly instead of $H_{t}^{N,N}$. More precisely, using the BBGKY hierarchy, the idea is to obtain a system of differential inequality of the form
\begin{equation} \frac{d}{dt} H_{t}^{k} \leq \frac{k (k-1)^{2}}{(N-1)^{2}} M+ \gamma k \bigl(H_{t}^{k+1}-H_{t}^{k}\bigl),
 \label{system 1} \end{equation}
and then to conclude relying on estimates on nested integrals.

Several results have since been obtained in this direction. Together with Le Flem, Lacker \cite{uniform in time} was able to obtain, under some more hypotheses, uniform in time results of order $\mathcal{O}\bigl(\frac{k^{2}}{N^{2}}\bigl)$. Their proof relies on the use of log-Sobolev inequalities to obtain an additional term of of the type $-c H_{t}^{k}$ in the right-hand side of \eqref{system 1}. Also, Lacker, Yeung and Zhou \cite{Lacker non exchangeable} were able to obtain sharp propagation of chaos results for systems of particles with weighted interactions (for instance, when the interaction is given by an interaction matrix). They compared the law of the system of particles to the law of n independent particles having the so called independent projection law, considering the BBGKY hierarchy in the case of non exchangeable diffusions. Moreover Hess-Childs and Rowan \cite{Hess-childs Rowan} obtained optimal rates for the $\chi^{2}$ divergence instead of the relative entropy, and they obtained also sharp convergence rates results for higher order corrections of the mean-field limit. 

Concerning singular interaction kernels Wang \cite{Wang} was recently able to obtain sharp rates for divergence free and $W^{-1,\infty}$ kernels (non optimal convergence for this case was obtained in the breakthrough paper \cite{Jabin wang}). A main step in the proof of this result is to consider a more general system of differential inequalities than \eqref{system 1}, of the form
\begin{equation}\label{system ineq Wang}
\frac{d}{dt} H_{t}^{k} \leq -c_{1} I_{t}^{k}+c_{2} I_{t}^{k+1} \mathbbm{1}_{k<N}+M_{1} H_{t}^{k} +M_{2} k \bigl(H_{t}^{k+1}-H_{t}^{k} \bigl) \mathbbm{1}_{k<N}+M_{3} e^{M_{3}t} \frac{k^{\beta}}{N^{2}},
\end{equation}
where $c_{1} \geq c_{2}\geq 0$ and $I_{t}^{k}:=\sum_{i=1}^{k} \int \fk \lVert \nabla_{v_{i}} \log \frac{\fk}{\ftk} \rVert^{2}$ is the Fisher information. With similar methods Wang was also able to obtain sharp rate for the $\chi^{2}$ divergence. Note moreover that, relying also on the BBGKY hierarchy but with a different approach (involving uniform in N estimates on weighted $L^{p}$ norms of the marginals), Bresch, Jabin and Soler \cite{Poisson Fokker Planck} were able to obtain convergence to the Vlasov-Poisson-Fokker-Planck equation (in dimension 2 and with a partial result for dimension 3) for short times.

\medskip

To the best of our knowledge, the result of this article provides for the first time an optimal rate of convergence for particles systems with a non constant diffusion coefficient. A central point of our proof is to obtain a system of inequalities of the form of \eqref{system ineq Wang}, controlling the additional terms appearing in the BBGKY hierarchy (the terms $K_{1}$, $K_{3}$ and $K_{4}$ of the proof). Our result is valid when the diffusion coefficient $a_2$ is sufficiently close to a constant. More precisely we prove the following Theorem. 


\begin{theorem}\label{th main}
Let $\bigl(V_{t}^{i,N} \bigl)_{i=1,..,N}$ solving (\ref{(1)}).
Let us assume that $b$, $\nabla \cdot a_{1}$, $a_{2}$ and $\nabla \cdot a_{2}$ are bounded. Suppose also the following: 
\begin{enumerate}
\label{i} \item[\emph{(i)}] {\sl Chaos at $t=0$:} there exists $C>0$ such that $H \bigl(\mu_{0}^{k,N}|\mu_{0}^{\otimes k} \bigl) \leq C \frac{k^{2}}{N^{2}}$.
\label{ii} \item[\emph{(ii)}] {\sl Exchangeability:} $\bigl( V_{0}^{1,N},...,V_{0}^{N,N} \bigl)$ are exchangeable.
\label{iii}\item[\emph{(iii)}]  {\sl Uniform ellipticity}: \eqref{condition a} is satisfied, i.e., there exists $\lambda_1>0$ such that
$$ \forall (v_{1},v_{2}) \in ({T}^{d})^{2},\,  a_{1}(v_{1})+a_{2}(v_{1}-v_{2}) \geq \lambda_{1} \text{Id}. $$
\label{iv} \item[\emph{(iv)}] {\sl Small dependency on the empirical measure for the diffusion coefficient }: $$\eta:=\underset{(x,z)\in (\mathbb{T}^{d})^{2}}{\sup} \lVert a_{2}(x)-a_{2}(z) \rVert_{2} < \lambda_{1}.$$
\end{enumerate} 
Then, for all $t>0$ there exists a constant $M_t$ independent of $N$ and $k$ such that
$$H_{t}^{k} \leq M_t \frac{k^{2}}{N^{2}}.$$
\end{theorem}
\textbf{{Remarks}}
\begin{itemize}
 \item Similarly as in \cite{Wang} one can prove a similar bound for the $\chi^{2}$ divergence instead of the entropy by adapting the proof of Theorem \ref{th main} (under very similar assumptions to (i), (ii), (iii) and (iv)). One can then obtain a system of the form
$$\frac{dD_{t}^{k}}{dt} \leq -c_{1} E_{t}^{k}+c_{2} E_{t}^{k+1} \mathbbm{1}_{k<N}+M_{1}D_{t}^{k+1} \mathbbm{1}_{k<N} +M_{2} \frac{k^{2}}{N^{2}},$$
where $D_{t}^{k}:=D(\fk|\ftk)$ is the $\chi^{2}$ divergence and the energy is defined as
$$E_{t}^{k}:=\sum_{i=1}^{k} \int \mu_{t}^{k,N} \lVert \nabla_{v_{i}} \frac{\mu_{t}^{k,N}}{\mu_{t}^{\otimes k}} \rVert^{2}.$$
Relying on Proposition 6 of \cite{Wang}, one obtains a bound of the form
$$D_{t}^{k} \leq \frac{M e^{M k}}{(T^{*}-t)^{3} N^{2}},$$
where M or $T^{*}$ do not depend on $N,k$ or $t$.
\item It is possible to relax the assumption on $b$ and consider a function of the form $b_{1}+b_{2}$, where $b_{1}=\nabla \cdot V \in W^{-1, \infty}(\mathbb{T}^{d},\mathbb{R}^{d})$ satisfies $\nabla \cdot b_{1}=0$, and $b_{2}$ is bounded, very much like in \cite{Wang}. By using Wang's estimates to control the terms involving $b_{1}$, it is possible to show that $H_{t}^{k}=\mathcal{O}\bigl(\frac{k^{2}}{N^{2}} \bigl)$ by controlling a system of the form \eqref{system ineq Wang}, where $c_{1}=\lambda_{1}-\sum_{i=1}^{4} \epsilon_{i}-\frac{\eta}{2 \epsilon}$ and $c_{2}= 2\eta \epsilon+\frac{(1+\epsilon_{5}) \lVert V \rVert^{2}_{\infty}}{4 \epsilon_{4}}$ for some $\epsilon_{i}, i \in \llbracket 1,5 \rrbracket$. The quantities $\epsilon_{1}, \epsilon_{2}, \epsilon_{3}, \epsilon_{5}$ can be chosen arbitrarily small (but not zero) and will only affect the value of M. In order to apply Wang's Proposition 5, we need to find $\epsilon_{4}$ such that $c_{1} > c_{2}$. By standard analysis, it is possible if and only if $\lambda_{1} > 2 \eta+ \lVert V \rVert_{\infty}$. This gives a generalisation of condition (iv).
\end{itemize}
\section{Proof of theorem \ref{th main}}
\subsection{Computation of the relative entropy and derivation of the system of differential inequalities}
The proof relies on the Fokker Planck equations satisfied by $\mu_{t}^{k,N}$ and $\mu_{t}^{\otimes k}$. The equation satisfied by $\mu_{t}^{k,N}$ is also known as the BBGKY hierarchy, because it gives $\partial_{t} \mu_{t}^{k,N}$ as a functional of $\mu_{t}^{k,N}$ and $\mu_{t}^{k+1,N}$. We defined $\hat{b}:=b-\nabla \cdot a_{2}$ to simplify the expressions. We have
\begin{align}
\nonumber \partial_{t} \mu_{t}^{k,N} = &\frac{1}{N}  \sum_{i=1}^{k} \sum_{j=1}^{k} \Bigl( -\nabla_{v_{i}} \cdot \left[ \hat{b}(v_{i}-v_{j}) \mu_{t}^{k,N} \right] + \nabla_{v_{i}} \cdot \left[ \bigl(a_{1}(v_{i})+a_{2}(v_{i}-v_{j}) \bigl)\nabla_{v_{i}} \mu_{t}^{k,N} \right] \Bigl)
\\ &
\nonumber+\sum_{i=1}^{k} \frac{N-k}{N} \int \nabla_{v_{i}} \cdot \left[ -\hat{b}(v_{i}-v_{k+1})  \mu_{t}^{k+1,N}\right] dv_{k+1}
\\ &
+\sum_{i=1}^{k} \frac{N-k}{N} \int \nabla_{v_{i}} \cdot \left[\bigl(a_{1}(v_{i})+a_{2}(v_{i}-v_{k+1}) \bigl )\nabla_{v_{i}} \mu_{t}^{k+1,N} \right] dv_{k+1}
\\ & \nonumber
+\sum_{i=1}^{k} \nabla_{v_{i}} \cdot \left[\nabla_{v_{i}} \cdot \bigl(a_{1}(v_{i})\bigl)
\mu_{t}^{k,N} \right],  
\end{align}
and
\begin{align}
\nonumber \partial_{t} \mu_{t}^{\otimes k}= \sum_{i=1}^{k} & \nabla_{v_{i}} \cdot \left[\bigl(a_{1}+a_{2}\ast \mu_{t}(v_{i}) \bigl) \nabla_{v_{i}} \mu_{t}^{\otimes k}\right]-\nabla_{v_{i}} \cdot \left[\bigl(\hat{b} \ast \mu_{t}(v_{i})\bigl) \mu_{t}^{\otimes k} \right]
\\ &
+\sum_{i=1}^{k} \nabla_{v_{i}} \cdot \left[ \nabla_{v_{i}} \cdot \bigl(a_{1}(v_{i}) \bigl) \mu_{t}^{\otimes k} \right]
\label{Fokker Planck limite}.
\end{align}
We can now compute the relative entropy:
\begin{align}
\nonumber \frac{d}{dt} H \bigl (\mu_{t}^{k,N}|\mu^{\otimes k }_{t} \bigl) & =\frac{d}{dt} \int \mu_{t}^{k,N} \log \Bigl(\frac{\mu_{t}^{k,N}}{\mu_{t}^{\otimes k }} \Bigl) 
\\ & \nonumber =\int \partial_{t} \mu_{t}^{k,N} \bigl(\log(\mu_{t}^{k,N})+1 \bigl)-\partial_{t} \mu_{t}^{k,N} \log(\mu_{t}^{\otimes k })-\partial_{t} \log(\mu_{t}^{\otimes k }) \mu_{t}^{k,N} \\ & =\mathcal{A}-\mathcal{B}-\mathcal{C}
\label{4}.
\end{align}
We have:
$$\mathcal{A}=\int \partial_{t} \mu_{t}^{k,N}+\partial_{t} \mu_{t}^{k,N} \log(\mu_{t}^{k,N})=\frac{d}{dt} \int \mu_{t}^{k,N}+\int \partial_{t} \mu_{t}^{k,N} \log(\mu_{t}^{k,N})=\int \partial_{t} \mu_{t}^{k,N} \log(\mu_{t}^{k,N}).$$
Therefore, by using the Fokker-Planck equations on $\mu_{t}^{\otimes k}$ and $\mu_{t}^{k,N}$, an integration by parts and the fact $\nabla \cdot a_{1}(v_{i})$ is a function of $v_{1},...,v_{k}$:
\begin{align}
\nonumber & \mathcal{A}-\mathcal{B} =-\sum_{i=1}^{k} \int \nabla_{v_{i}} \cdot \bigl(a_{1}(v_{i})\bigl) \fk \cdot \nabla_{v_{i}}\log \left( \frac{\fk}{\ftk} \right) 
 \\ & \nonumber+\int \sum_{i=1}^{k} \frac{1}{N} \left[-\sum_{j=1}^{k} \bigl( a_{1}(v_{i})+a_{2}(v_{i}-v_{j}) \bigl) \nabla_{v_{i}} \fk +\sum_{j=1}^{k} \hat{b}(v_{i}-v_{j}) \fk \right] \cdot  \nabla_{v_{i}}\log \left( \frac{\fk}{\ftk} \right) 
 \\ & +\sum_{i=1}^{k} \frac{N-k}{N} \int  -\bigl( a_{1}(v_{i})+a_{2}(v_{i}-v_{k+1}) \bigl) \nabla_{v_{i}} \mu_{t}^{k+1,N} \cdot \nabla_{v_{i}}\log \left( \frac{\fk}{\ftk} \right) \label{5}
\\ & \nonumber
+\sum_{i=1}^{k} \frac{N-k}{N} \int \hat{b}(v_{i}-v_{k+1}) \mu_{t}^{k+1,N} \cdot \nabla_{v_{i}}\log \left( \frac{\fk}{\ftk} \right).
\end{align}
We have, by integration by parts and the identity
$\nabla_{v_{i}} \Bigl(\frac{\fk}{\ftk}\Bigl)=\frac{\fk}{\ftk} \nabla_{v_{i}} \log \frac{\fk}{\ftk}$:
\begin{align}
 \nonumber \mathcal{C}=\sum_{i=1}^{k}& \int  - \bigl(a_{1}(v_{i})+a_{2}\ast \mu_{t}(v_{i}) \bigl) \nabla_{v_{i}} \log(\ftk)\cdot \fk \nabla_{v_{i}} \log \frac{\fk}{\ftk}
\\ & +\sum_{i=1}^{k} \int \bigl(\hat{b} \ast \mu_{t}(v_{i})- \nabla_{v_{i}} \cdot a_{1}(v_{i})\bigl) \cdot \fk \nabla_{v_{i}} \log \frac{\fk}{\ftk} \label{6}.
\end{align}
Combining (\ref{4}), (\ref{5}) and (\ref{6}), rearranging the terms and using the identities 
$$\nabla_{v_{i}} \ftk=\ftk \nabla_{v_{i}} \log \ftk \quad\mbox{ and }\quad \nabla_{v_{i}} \log \fk=\nabla_{v_{i}} \log \ftk+\nabla_{v_{i}} \log \frac{\fk}{\ftk},$$
we deduce
\begin{align*}
\frac{d}{dt} H_{t}^{k} =&\sum_{i=1}^{k}  \int \fk \left[\frac{1}{N} \sum_{j=1} ^{k} \bigl(a_{1}(v_{i})+a_{2}(v_{i}-v_{j}) \bigl)\right]\nabla_{v_{i}} \log (\ftk) \cdot \nabla_{v_{i}} \log \frac{\fk}{\ftk}  \\ &- \sum_{i=1}^{k} \int \fk \bigl(a_{1}(v_{i})+a_{2} \ast \mu_{t}(v_{i}) \bigl) \nabla_{v_{i}} \log (\ftk) \cdot \nabla_{v_{i}} \log \frac{\fk}{\ftk} 
\\ & 
-\frac{1}{N} \sum_{i=1}^{k} \int \fk \left[\sum_{j=1} ^{k} \bigl(a_{1}(v_{i})+a_{2}(v_{i}-v_{j}) \bigl) \right] \nabla_{v_{i}} \log \frac{\fk}{\ftk} \cdot \nabla_{v_{i}} \log \frac{\fk}{\ftk} 
\\ & 
+ \sum_{i=1} ^{k} \int \fk  \left[ \frac{1}{N}\sum_{j=1} ^{k}  \bigl(\hat{b}(v_{i}-v_{j})\bigl)-\hat{b} \ast \mu_{t}(v_{i}) \right] \cdot \nabla_{v_{i}} \log \frac{\fk}{\ftk} 
\\ &
-\frac{N-k}{N}\sum_{i=1}^{k} \int \bigl(a_{1}(v_{i})+a_{2}(v_{i}-v_{k+1})\bigl)\nabla_{v_{i}} \mu_{t}^{k+1,N} \cdot \nabla_{v_{i}} \log \Bigl(\frac{\fk}{\ftk} \Bigl)
\\ &
-\frac{N-k}{N} \sum_{i=1}^{k} \int \hat{b}(v_{i}-v_{k+1}) \mu_{t}^{k+1,N} \cdot \nabla_{v_{i}} \log \Bigl(\frac{\fk}{\ftk} \Bigl)
\end{align*}
Because of the assumption on $a_{1}$ and $a_{2}$ we have:
\begin{align*}
 -\sum_{i=1}^{k} \int \fk & \left[\sum_{j=1} ^{k} \frac{1}{N} \bigl(a_{1}(v_{i})+a_{2}(v_{i}-v_{j}) \bigl) \right] \nabla_{v_{i}} \log \frac{\fk}{\ftk} \cdot \nabla_{v_{i}} \log \frac{\fk}{\ftk} 
\\ &
\leq -\lambda_{1} \frac{k}{N} \sum_{i=1}^{k} \int \fk \Big\| \nabla_{v_{i}} \log \frac{\fk}{\ftk} 
 \Big\| ^{2} .
\end{align*}
By using the fact that $a_{1}(v_{i})+a_{2} \ast \mu_{t}(v_{i})=\frac{k}{N} \bigl(a_{1}(v_{i})+a_{2} \ast \mu_{t}(v_{i}) \bigl)+\frac{N-k}{N} \bigl(a_{1}(v_{i})+a_{2} \ast \mu_{t}(v_{i}) \bigl)$, rearranging the terms and by using the previous inequality:
\begin{align}
\nonumber \frac{d}{dt} H_{t}^{k} & \leq -\sum_{i=1}^{k} \int \fk \Bigr[ \frac{1}{N} \sum_{j=1} ^{k} \bigl(a_{2} \ast \mu_{t}(v_{i})-a_{2}(v_{i}-v_{j}) \bigl) \Bigr]  \nabla_{v_{i}} \log (\ftk) \cdot \nabla_{v_{i}} \log \frac{\fk}{\ftk} 
\\ & \nonumber
\ \ \ \ +\sum_{i=1} ^{k} \int \fk \left[\sum_{j=1} ^{k} \frac{1}{N} \hat{b}(v_{i}-v_{j})-\frac{k}{N} \hat{b}\ast \mu_{t}(v_{i}) \right] \cdot \nabla_{v_{i}} \log \frac{\fk}{\ftk} 
\\ & \nonumber
\ \ \ \ -\lambda_{1} \frac{k}{N} \sum_{i=1}^{k} \int \fk \Big\| \nabla_{v_{i}} \log \frac{\fk}{\ftk} \Big\|^{2}-\sum_{i=1}^{k} \frac{N-k}{N} \int \fk \hat{b} \ast \mu_{t}(v_{i}) \cdot \nabla_{v_{i}} \log \frac{\fk}{\ftk}
\\ & \nonumber
\ \ \ \ -\sum_{i=1}^{k} \frac{N-k}{N} \int \bigl(a_{1}(v_{i})+a_{2}(v_{i}-v_{k+1}) \bigl) \nabla_{v_{i}} \mu_{t}^{k+1,N}  \cdot \nabla_{v_{i}} \log \frac{\fk}{\ftk}
\\ & \nonumber
\ \ \ \ -\sum_{i=1}^{k} \frac{N-k}{N} \int \hat{b}(v_{i}-v_{k+1}) \mu_{t}^{k+1,N} \cdot \nabla_{v_{i}} \log \frac{\fk}{\ftk}
\\ & \nonumber
\ \ \ \ +\sum_{i=1}^{k} \frac{N-k}{N} \int \fk \bigl(a_{1}(v_{i})+a_{2} \ast \mu_{t}(v_{i}) \bigl) \nabla_{v_{i}} \log (\ftk) \cdot \nabla_{v_{i}} \log \frac{\fk}{\ftk} 
\\ &
\leq  I+J -\lambda_{1} \frac{k}{N} \sum_{i=1}^{k} \int \fk \Big\| \nabla_{v_{i}} \log \frac{\fk}{\ftk} \Big\|^{2}+K \label{7}.
\end{align}
So, using the classical inequality $x \cdot y \leq N\epsilon_{1} \lVert x \rVert ^{2}+\frac{1}{4\epsilon_{1} N} \lVert y \rVert ^{2}$ we have, for all $\epsilon_{1}>0$,
\begin{align}
 \nonumber I & = \sum_{i=1}^{k} \frac{1}{N} \int \fk \left[\sum_{j=1}^{k} \bigl(a_{2}(v_{i}-v_{j})-a_{2} \ast \mu_{t}(v_{i}) \bigl) \right] \nabla_{v_{i}} \log(\ftk) \cdot \nabla_{v_{i}} \log \frac{\fk}{\ftk}
\\ & \nonumber \leq \sum_{i=1}^{k} \epsilon_{1} \int \fk \Big\| \nabla_{v_{i}} \log \frac{\fk}{\ftk} \Big\|^{2}+\frac{1}{4\epsilon_{1} N^{2}} C_t^{2} \int \mu_{t}^{k,N} \big\| \sum_{j=1}^{k} \bigl(a_{2}(v_{i}-v_{j})-a_{2}\ast \mu_{t}(v_{i}) \bigl) \big\|^{2} 
\\ &
\leq \sum_{i=1}^{k} \epsilon_{1} \int \fk \Big\| \nabla_{v_{i}} \log \frac{\fk}{\ftk} \Big\|^{2} + C_t^{2}  \frac{C k^{3}}{\epsilon_{1} N^{2}} \label{8}.
\end{align}
 Note that we used the fact that $\sup_{s\le t}\|\nabla \log \mu_s\|\le C_t$ to get the second inequality. We refer to \cite[Prop. 3.1]{BGP10} where this statement is proved in the case of a constant diffusion case, but the proof can readily be extended to the uniformly elliptic bounded diffusion case. $C_{t}$ is fixed but the value of $C$ may change from line to line and it will be independent of $N$, $k$ and t. By using the same type of estimates, we can obtain a very similar bound for $J$:
\begin{align} 
J&=\sum_{i=1} ^{k} \int \fk \left[\sum_{j=1} ^{k} \frac{1}{N} \hat{b}(v_{i}-v_{j})-\frac{k}{N} \hat{b}\ast \mu_{t}(v_{i}) \right] \cdot \nabla_{v_{i}} \log \frac{\fk}{\ftk} \\
&\leq \sum_{i=1}^{k} \epsilon_{2} \int \fk \Big\| \nabla_{v_{i}} \log \frac{\fk}{\ftk} \Big\|^{2} + C \frac{k^{3}}{\epsilon_{2} N^{2}} \label{9}. \end{align}
Let us recall the expression of $K$:
\begin{align*}K =& -\sum_{i=1}^{k} \frac{N-k}{N} \int \bigl(a_{1}(v_{i})+a_{2}(v_{i}-v_{k+1}) \bigl)\nabla_{v_{i}} \mu_{t}^{k+1,N} \cdot \nabla_{v_{i}} \log \frac{\fk}{\ftk}
\\ &
+\sum_{i=1}^{k} \frac{N-k}{N} \int \fk \bigl(a_{1}(v_{i})+a_{2} \ast \mu_{t}(v_{i}) \bigl) \nabla_{v_{i}} \log (\ftk) \cdot \nabla_{v_{i}} \log \frac{\fk}{\ftk} 
\\ &-\sum_{i=1}^{k} \frac{N-k}{N} \int \Bigl(\fk \hat{b} \ast \mu_{t}(v_{i}) \cdot \nabla_{v_{i}} \log \frac{\fk}{\ftk}+ \hat{b}(v_{i}-v_{k+1}) \mu_{t}^{k+1,N} \cdot \nabla_{v_{i}} \log \frac{\fk}{\ftk} \Bigl).
\end{align*}
We can simplify the terms involving $\hat{b}$ because of the following remark:
$$\int \fk \hat{b} \ast \mu_{t}(v_{i}) \cdot \nabla_{v_{i}} \log \frac{\fk}{\ftk}=\int\mu_{t}^{k+1,N} \hat{b} \ast \mu_{t}(v_{i}) \cdot \nabla_{v_{i}} \log \frac{\fk}{\ftk}.$$
This is valid because $\mu_{t}^{k+1,N}$ is integrated against a function of $v_{1},...,v_{k}$ only. We can also simplify the terms involving $a_{1}$ and $a_{2}$ as follows:
\begin{align*}
& \int \fk \bigl(a_{1}(v_{i})+a_{2} \ast \mu_{t}(v_{i}) \bigl) \nabla_{v_{i}} \log (\ftk) \cdot \nabla_{v_{i}} \log \frac{\fk}{\ftk} \\ & = \int \fk \bigl(a_{1}(v_{i})+a_{2} \ast \mu_{t}(v_{i}) \bigl) \nabla_{v_{i}} \log (\fk) \cdot \nabla_{v_{i}} \log \frac{\fk}{\ftk}
\\ &
\ \ 
-\int \fk \bigl(a_{1}(v_{i})+a_{2} \ast \mu_{t}(v_{i})\bigl) \nabla_{v_{i}} \log (\frac{\fk}{\ftk}) \cdot \nabla_{v_{i}} \log \frac{\fk}{\ftk}
\\ &
\leq \int \bigl(a_{1}(v_{i})+a_{2} \ast \mu_{t}(v_{i}) \bigl) \nabla_{v_{i}} \fk \cdot \nabla_{v_{i}} \log \frac{\fk}{\ftk}
-\lambda_{1} \int \fk \Big\| \nabla_{v_{i}} \log \frac{\fk}{\ftk} \Big\|^{2}.
\end{align*}
Let
$$\mathcal{D}=\int \bigl(a_{1}(v_{i})+a_{2} \ast \mu_{t}(v_{i}) \bigl)\nabla_{v_{i}} \fk \cdot \nabla_{v_{i}} \log \frac{\fk}{\ftk}.$$
Using a double integration by parts, we obtain
\begin{align*}
\mathcal{D} & = \int \bigl(a_{1}(v_{i})+a_{2} \ast \mu_{t}(v_{i}) \bigl) \nabla_{v_{i}} \log \frac{\fk}{\ftk}
\cdot \nabla_{v_{i}} \fk
\\ &
=-\int \fk \nabla_{v_{i}} \cdot \left[ \bigl(a_{1}(v_{i})+a_{2} \ast \mu_{t}(v_{i}) \bigl) \nabla_{v_{i}} \log \frac{\fk}{\ftk} \right]
\\ & =-\int \mu_{t}^{k+1,N} \nabla_{v_{i}} \cdot \left[ \bigl(a_{1}(v_{i})+a_{2} \ast \mu_{t}(v_{i}) \bigl) \nabla_{v_{i}} \log \frac{\fk}{\ftk} \right]
\\ &
=\int \nabla_{
v_{i}} \mu_{t}^{k+1,N} \bigl(a_{1}(v_{i})+a_{2} \ast \mu_{t}(v_{i}) \bigl) \nabla_{v_{i}} \log \frac{\fk}{\ftk}.
\end{align*}
\vspace{-0,2cm}
Combining these estimates and simplifying the expressions, we get
\begin{align}
\nonumber K & \leq -\sum_{i=1}^{k} \frac{N-k}{N} \lambda_{1} \int \fk \Big\| \nabla_{v_{i}} \log \frac{\fk}{\ftk} \Big\|^{2}
\\ & \nonumber \ \ \ \ +\sum_{i=1}^{k} \frac{N-k}{N} \int \bigl(a_{2} \ast \mu_{t}(v_{i})-a_{2}(v_{i}-v_{k+1}) \bigl) \nabla_{v_{i}} \mu_{t}^{k+1,N} \cdot \nabla_{v_{i}} \log \frac{\fk}{\ftk}
\\ & \nonumber
\ \ \ \ +\sum_{i=1}^{k} \frac{N-k}{N} \int \bigl(\hat{b}(v_{i}-v_{k+1})-\hat{b} \ast \mu_{t}(v_{i}) \bigl)\mu_{t}^{k+1,N} \cdot \nabla_{v_{i}} \log \frac{\fk}{\ftk}
\\ &
\leq -\sum_{i=1}^{k} \frac{N-k}{N} \lambda_{1} \int \fk \Big\| \nabla_{v_{i}} \log \frac{\fk}{\ftk} \Big\|^{2}+K_{1}+K_{2}\label{10}.
\end{align}
We are left with estimating the terms $K_{1}$ and $K_{2}$.

On one hand, integrating with respect to $v_{1},...,v_{k}$ first then by respect to $v_{k+1}$, and by multiplying and dividing by $\mu_{t}^{k,N}(v_{1},...,v_{k})$, we have the following expression for $K_{2}$:
\begin{multline*}
K_{2} = \sum_{i=1}^{k} \frac{N-k}{N} \int_{v_{1},...,v_{k}} \fk\nabla_{v_{i}} \log \frac{\fk}{\ftk} 
\\ \cdot \Bigr \langle \hat{b}(v_{i}-.),v_{k+1} \mapsto \mu_{t}^{k+1|k,N}(v_{k+1}|v_{1},...,v_{k})-\mu_{t} \Bigr \rangle,
\end{multline*}
where $v_{k+1} \mapsto \mu_{t}^{k+1|k,N}(v_{k+1}|v_{1},...,v_{k})$ is the conditonal density of $V_{1},...,V_{k+1}$ with respect to $(V_{1}=v_{1},...,V_{k}=v_{k})$.
Pinsker inequality yields, for all $v_{1},...,v_{k}$:
$$\Big\| \langle \hat{b}(v_{i}-.),v_{k+1} \mapsto \mu_{t}^{k+1|k,N}(v_{k+1}|v_{1},...,v_{k})-\mu_{t} \rangle \Big\| \leq \lVert \hat{b} \rVert_{\infty} \sqrt{H\bigl(\mu_{t}^{k+1|k,N}(v_{k+1}|v_{1},...,v_{k})|\mu_{t} \bigl)}.$$
Finally, by using $x\cdot y \leq \epsilon \lVert x \rVert ^{2}+\frac{1}{4\epsilon} \lVert y \rVert ^{2}$, we get:
\begin{align*}
K_{2} & \leq \sum_{i=1}^{k} \frac{N-k}{N} \int \epsilon \fk \Big\| \nabla_{v_{i}} \log \frac{\fk}{\ftk} \Big\|^{2}
\\ & 
\ \ \ \ +\sum_{i=1}^{k} \frac{N-k}{4 \epsilon N} \int \fk \Big\| \Bigr \langle \hat{b}(v_{i}-.),v_{k+1} \mapsto \mu_{t}^{k+1|k,N}(v_{k+1}|v_{1},...,v_{k})-\mu_{t} \Bigr \rangle \Big\|^{2}
\\ &
\leq \sum_{i=1}^{k} \frac{N-k}{N} \int \epsilon \fk \Big\| \nabla_{v_{i}} \log \frac{\fk}{\ftk} \Big\|^{2}
\\ &
\ \ \ \ +\sum_{i=1}^{k} \frac{N-k}{4 \epsilon N} \int \fk \lVert \hat{b} \rVert_{\infty}^{2} H\bigl(\mu_{t}^{k+1|k,N}(v_{k+1}|v_{1},...,v_{k})|\mu_{t} \bigl).
\end{align*}
By a towering property of the relative entropy
$$\int \fk H\bigl(\mu_{t}^{k+1|k,N}(v_{k+1}|v_{1},...,v_{k})|\mu_{t} \bigl)=H_{t}^{k+1}-H_{t}^{k},$$ 
which can easily be checked by using the definition of relative entropy, we get:
\begin{align}
K_{2} & \leq \sum_{i=1}^{k} \frac{N-k}{N} \int \epsilon \fk \Big\| \nabla_{v_{i}} \log \frac{\fk}{\ftk} \Big\|^{2}
+ k\frac{C(N-k)}{4\epsilon N}\bigl(H_{t}^{k+1}-H_{t}^{k} \bigl). \label{11}
\end{align}

On the other hand, we can rewrite $K_{1}=K_{3}+K_{4}$, with 
\begin{align*}
& K_{3} =\sum_{i=1}^{k} \frac{N-k}{N} \int \mu_{t}^{k+1,N} \bigl(a_{2} \ast \mu_{t}(v_{i})-a_{2}(v_{i}-v_{k+1}) \bigl) \nabla_{v_{i}} \log (\mu_{t}^{ \otimes k+1}) \cdot \nabla_{v_{i}} \log \frac{\fk}{\ftk}
\\ &
K_{4}=\sum_{i=1}^{k} \frac{N-k}{N} \int \mu_{t}^{k+1,N} \bigl(a_{2} \ast \mu_{t}(v_{i})-a_{2}(v_{i}-v_{k+1}) \bigl) \nabla_{v_{i}} \log \Bigl(\frac{\mu_{t}^{ k+1,N}}{\mu_{t}^{\otimes k+1}} \Bigl) \cdot \nabla_{v_{i}} \log \frac{\fk}{\ftk}.
\end{align*}
The key observation is the following:
$$\nabla_{v_{i}} \log \mu_{t}^{\otimes k+1}=\nabla_{v_{i}} \left[\sum_{j=1}^{k+1} \log \bigl(\mu_{t}(v_{j}) \bigl) \right]=\nabla_{v_{i}} \left[ \log \bigl(\mu_{t}(v_{i}) \bigl) \right]=\nabla_{v_{i}} \log \mu_{t}^{\otimes k}.
$$
Therefore, by integrating with respect to $v_{1},...,v_{k}$ first:
\begin{align*}
K_{3} & =\sum_{i=1}^{k} \frac{N-k}{N} \int_{v_{1},...,v_{k}} \fk \nabla_{v_{i}} \log \frac{\fk}{\ftk} \cdot \Bigr \langle a_{2}(v_{i}-.),\frac{\mu_{t}^{k+1,N}}{\fk}-\mu_{t} \Bigr \rangle \nabla_{v_{i}} \log \ftk .
\end{align*}
We can now use a very similar argument than what we used to control $K_{2}$, by using the fact that $\nabla_{v_{i}} \log \mu_{t}^{\otimes k}$ to get:
\begin{align}
K_{3}\leq \sum_{i=1}^{k} \frac{N-k}{N} \epsilon \int \fk \Big\| \nabla_{v_{i}} \log \frac{\fk}{\ftk} \Big\|^{2}+k C_{t}^{2} \frac{C(N-k)}{4 \epsilon N} \bigl(H_{t}^{k+1}-H_{t}^{k} \bigl)  \label{12}.
\end{align}
We can now control the term $K_{4}$, for which we will rely on condition (iv) of Theorem~\ref{th main}. By the Cauchy-Schwarz inequality and by definition of $a_{2} \ast \mu_{t}(v_{i})$, we have:
\begin{multline*}
   K_{4} \leq \sum_{i=1}^{k} \frac{N-k}{N} \int \mu_{t}^{k+1,N}  \Big\| \Bigl( \int \bigl(a_{2}(v_{i}-v) -a_{2}(v_{i}-v_{k+1}) \bigl) \mu_{t}(v) \Bigl) \nabla_{v_{i}} \log \frac{\mu_{t}^{ k+1,N}}{\mu_{t}^{\otimes k+1}} \Big\| 
\\
\times \Big\| \nabla_{v_{i}} \log \frac{\fk}{\ftk} \Big\| .
\end{multline*}
Therefore, by the inequality $x \cdot y \leq \epsilon \lVert x \rVert^{2}+ \frac{1}{4 \epsilon} \lVert y \rVert^{2}$:
\begin{align*}
K_{4} &\leq \sum_{i=1}^{k} \frac{N-k}{N} \int \mu_{t}^{k+1,N}    \eta \Big\| \nabla_{v_{i}} \log \frac{\mu_{t}^{ k+1,N}}{\mu_{t}^{\otimes k+1}} \Big\| \Big\| \nabla_{v_{i}} \log \frac{\fk}{\ftk} \Big\|
\\ &
\leq \sum_{i=1}^{k}  \eta \epsilon \frac{N-k}{N} \int \mu_{t}^{k+1,N} \Big\| \nabla_{v_{i}} \log \frac{\mu_{t}^{ k+1,N}}{\mu_{t}^{\otimes k+1}} \Big\|^{2}+\sum_{i=1}^{k}  \eta  \frac{N-k}{4 \epsilon N} \int \mu_{t}^{k,N} \Big\| \nabla_{v_{i}} \log \frac{\mu_{t}^{ k,N}}{\mu_{t}^{\otimes k}} \Big\|^{2}.
\end{align*}
By combining (\ref{7}), (\ref{8}), (\ref{9}), (\ref{10}), (\ref{11}), (\ref{12}) we get:
\begin{align}
\nonumber \frac{d}{dt} H_{t}^{k}  \leq \bigl(-\lambda_{1} & +\epsilon_{1}+\epsilon_{2}+\epsilon_{3}+\frac{\eta}{4\epsilon} \bigl) \sum_{i=1}^{k} \int \fk \Big\| \nabla_{v_{i}} \log \frac{\fk}{\ftk} \Big\|^{2}
\\ & \nonumber +\sum_{i=1}^{k+1}  \eta \epsilon  \int \mu_{t}^{k+1,N} \Big\| \nabla_{v_{i}} \log \frac{\mu_{t}^{ k+1,N}}{\mu_{t}^{\otimes k+1}} \Big\|^{2}
+(C_t^{2}+1) C\frac{k^{3}}{N^{2}}
\\ &+(C_{t}^{2}+1) C k \bigl(H_{t}^{k+1}-H_{t}^{k} \bigl) \label{13}.
\end{align}
\subsection{Estimation of the system of inequalities}

Let us see how to use \eqref{13} to get Theorem~\ref{th main}.
We will mainly use the Proposition 5 from Wang \cite{Wang}, rewritten here:
\begin{proposition}\label{prop-wang}
Let $T>0$ and let $x_{.}^{k},y_{.}^{k}: [0,T] \to \mathbb{R}_{\geq 0}$ be $\mathcal{C}^{1}$ functions, for $k\in \llbracket 1,N \rrbracket$. Suppose that $x_{t}^{k+1} \geq x_{t}^{k} $ for all $k \in \llbracket 1,N-1 \rrbracket$. Suppose that there exists $\beta \geq 2$, real numbers $c_{1}>c_{2}\geq 0$ and $C_{0},M_{1},M_{2},M_{3} \geq 0$ such that, for all $t \in [0,T]$ and $k \in \llbracket 1,N \rrbracket$, we have
$$x_{0}^{k} \leq C_{0} \frac{k^{2}}{N^{2}}$$
$$\frac{dx_{t}^{k}}{dt} \leq -c_{1} y_{t}^{k}+c_{2} y_{t}^{k+1} \mathbbm{1}_{k<N}+M_{1} x_{t}^{k} +M_{2} k \bigl(x_{t}^{k+1}-x_{t}^{k} \bigl) \mathbbm{1}_{k<N}+M_{3} e^{M_{3}t} \frac{k^{\beta}}{N^{2}}.$$
Then, there exists $M>0$ depending only on $\beta,c_{1},c_{2}, C_{0}$ and $M_{i}, i \in \{1,2,3 \}$ such that, for all $t \in [0,T]$, we have
$$x_{t}^{k} \leq \frac{M e^{Mt} k^{\beta}}{N^{2}}.$$
\end{proposition}
Remark that (\ref{13}) is exactly of this form with $x_{t}^{k}=H_{t}^{k}$ and $y_{t}^{k}=I_{t}^{k}$ and with $M_{1}=0, \beta=3$. The inequality $H_{t}^{k+1} \geq H_{t}^{k}$ holds because of the towering property of the relative entropy. The value of our constants depends on t but we can use the fact that $t \mapsto C_{t}$ is bounded on $[0,T]$ for all $T>0$. In order to use the result, we need to ensure that $c_{1} > c_{2}$. By letting $\epsilon':=\epsilon_{1}+\epsilon_{2}+\epsilon_{3}$, it is equivalent to finding $\epsilon>0$ such that:
$$\lambda_{1}-\epsilon'-\frac{\eta}{4 \epsilon} >  \eta \epsilon.$$
It is possible if and only if $\eta$ satisfies (iii). 
We can now apply Wang's result to get the suboptimal estimate:
$$\forall T>0, \exists M(T) \in (0,+\infty), \forall t \in [0,T], H_{t}^{k} \leq \frac{M(T) e^{M(T)t} k^{3}}{N^{2}}.$$
We can improve this result by using it to prove better bounds for $I$ and $J$, as done in \cite{Lacker} or \cite{Wang}.
\newline 
Let us recall one of the previous bound we obtained for I:
\begin{align}
I \leq \sum_{i=1}^{k} \frac{\epsilon}{N} \int \fk \Big\| \nabla_{v_{i}} \log \frac{\fk}{\ftk} \Big\|^{2}+\frac{1}{4\epsilon N} C_{t}^{2} \int \mu_{t}^{k,N} \Big\| \sum_{j=1}^{k} \bigl(a_{2}(v_{i}-v_{j})-a_{2}\ast \mu_{t}(v_{i}) \bigl) \Big\|^{2}
 \label{14}.
\end{align}
We used the fact that $a_{2}$ is bounded to get the first bound, but we can develop the norm to get a better estimate. Indeed, we have:
\begin{align*}
& \int \mu_{t}^{k,N} \Big\| \sum_{j=1}^{k} \bigl(a_{2}(v_{i}-v_{j})-a_{2}\ast \mu_{t}(v_{i}) \bigl) \Big\|^{2}  \\ & \quad= \int \fk \sum_{j=1}^{k} \Big\| a_{2}(v_{i}-v_{j})-a_{2} \ast \mu_{t}(v_{i}) \Big\|^{2}
\\ &
\ \ \ \ \quad+\int \fk \sum_{j_{1} \neq j_{2}} \bigl( a_{2}(v_{i}-v_{j_{1}})-a_{2} \ast \mu_{t}(v_{i}) \bigl) \cdot \bigl(a_{2}(v_{i}-v_{j_{2}})-a_{2} \ast \mu_{t}(v_{i}) \bigl).
\end{align*}
Therefore:
\begin{align}
\nonumber & \int \mu_{t}^{k,N} \Big\| \sum_{j=1}^{k} \bigl(a_{2}(v_{i}-v_{j})-a_{2}\ast \mu_{t}(v_{i}) \bigl) \Big\|^{2} \\ & \quad\nonumber \leq \sum_{j_{1} \neq j_{2}} \int \fk \bigl( a_{2}(v_{i}-v_{j_{1}})-a_{2} \ast \mu_{t}(v_{i})\bigl) \cdot \bigl(a_{2}(v_{i}-v_{j_{2}})-a_{2} \ast \mu_{t}(v_{i}) \bigl)
+4 k \rVert a_{2} \lVert_{\infty}^{2} 
\\ &\quad
\leq k^{2} \int \mu_{t}^{3,N} \bigl( a_{2}(v_{1}-v_{2})-a_{2} \ast \mu_{t}(v_{1})) \bigl) \cdot \bigl(a_{2}(v_{1}-v_{3})-a_{2} \ast \mu_{t}(v_{1}) \bigl)
+4 k \rVert a_{2} \lVert_{\infty}^{2}
 \label{15}.
\end{align}
The last step is due to the exchangeability of the particles. Notice that we have:
\begin{align}
\nonumber \int \mu_{t}^{\otimes 3}& \bigl( a_{2}(v_{1}-v_{2})-a_{2} \ast \mu_{t}(v_{1}) \bigl) \cdot \bigl(a_{2}(v_{1}-v_{3})-a_{2} \ast \mu_{t}(v_{1}) \bigl) \\ \nonumber& = \int_{v_{1},v_{2}} \mu_{t}^{\otimes 2} \bigl(a_{2}(v_{1}-v_{2})-a_{2} \ast \mu_{t}(v_{1}) \bigl)  \cdot \left[ \int_{v_{3}} \bigl( a_{2}(v_{1}-v_{3})-a_{2} \ast \mu_{t}(v_{i}) \bigl) \mu_{t}(v_{3}) \right]\\
&=0
 \label{16},
\end{align}
as the term of integral with respect to $v_{3}$ is equal to zero.
 Therefore, combining (\ref{15}), (\ref{16}) and using the Pinsker inequality, we have:
\begin{align*}\int  & \mu_{t}^{k,N} \Big\| \sum_{j=1}^{k} \bigl(a_{2}(v_{i}-v_{j})-a_{2}\ast \mu_{t}(v_{i}) \bigl) \Big\|^{2} 
\\ &
\leq k^{2} \left[ \int \mu_{t}^{\otimes 3} \bigl( a_{2}(v_{1}-v_{2})-a_{2} \ast \mu_{t}(v_{1}) \bigl) \cdot \bigl(a_{2}(v_{1}-v_{3})-a_{2} \ast \mu_{t}(v_{1}) \bigl)  +C\sqrt{H\bigl(\mu_{t}^{3,N}|\mu_{t}^{\otimes 3}\bigl)} \right]
\\ &\quad+4k \lVert  a_{2} \rVert_{\infty}^{2}
\\ &
\leq k^{2} C\sqrt{H\bigl(\mu_{t}^{3,N}|\mu_{t}^{\otimes 3} \bigl)}+4k \lVert a_{2} \rVert_{\infty}^{2}.
\end{align*}
By combining this with (\ref{14}), we get:
\begin{align*}
I & \leq  \sum_{i=1}^{k} \epsilon_{1}\int \fk \Big\| \nabla_{v_{i}} \log \frac{\fk}{\ftk} \Big\|^{2}+\frac{C}{4 \epsilon_{1} N^{2}} C_{t}^{2} \left(k+k^{2} \sqrt{H \bigl(m^{3}|m^{\otimes 3} \bigl)}\right)
\\ & 
\leq \sum_{i=1}^{k} \epsilon_{1} \int \fk \Big\| \nabla_{v_{i}} \log \frac{\fk}{\ftk} \Big\|^{2} +C_{t}^{2}  \frac{Ck^{2}}{N^{2}}.
\end{align*}
\newline 
Using a similar argument for $J$ yields the following bound for $H_{t}^{k}$:
\begin{align*}
\frac{d}{dt} H_{t}^{k} & \leq -c_{1} \sum_{i=1}^{k} \int \fk \Big\| \nabla_{v_{i}} \log \frac{\fk}{\ftk} \Big\|^{2}+\sum_{i=1}^{k+1} c_{2}  \int \mu_{t}^{k+1,N} \Big\| \nabla_{v_{i}} \log \frac{\mu_{t}^{ k+1,N}}{\mu_{t}^{\otimes k+1}} \Big\|^{2}
\\ &
\ \ \ \ +(C_{t}^{2}+1) C\frac{k^{2}}{N^{2}}+(C_{t}^{2}+1) Ck \bigl(H_{t}^{k+1}-H_{t}^{k} \bigl). 
\end{align*}
Using Proposition~\ref{prop-wang} one more time leads us to the result and finishes the proof of our main result.

\medskip

\noindent{\bf Acknowledgments.}\\
This work has been (partially) supported by the Project CONVIVIALITY ANR-23-CE40-0003 and the Project PERISTOCH ANR–19–CE40–0023 of the French National Research Agency. AG has benefited from a government grant managed by the Agence Nationale de la
Recherche under the France 2030 investment plan "ANR-23-EXMA-0001”.

\bibliographystyle{amsplain}

\end{document}